\documentclass{ajour}

\usepackage{amssymb,amsmath}
\usepackage[mathscr]{eucal}
\DeclareMathOperator*{\Int}{\int}

\begin{document}

%
\finaltypesetting \journame{Journal of Combinatorial Theory,
Series A} \articlenumber{doi:10.1006/jcta.2001.3194}
\yearofpublication{2002} \volume{97} \cccline{0097-3165/02
\$35.00} \received{October 3, 2000}

\commline{Communicated by the Managing Editors}

\authorrunninghead{D.\ Bowman \& D.\ M.\ Bradley}
\titlerunninghead{Shuffles and Multiple Zeta Values}

\setcounter{page}{43} 




\title{The Algebra and Combinatorics of
Shuffles and Multiple Zeta Values}


%

%
\author{Douglas Bowman}
\affil{Department of Mathematics, Northern Illinois University,
DeKalb, Illinois} \email{bowman@math.niu.edu} \and

\author{David M. Bradley}
\affil{Department of Mathematics and Statistics, University of
Maine, Orono, Maine} \email{bradley@math.umaine.edu,
dbradley@e-math.ams.org}






\abstract{The algebraic and combinatorial theory of shuffles,
introduced by Chen and Ree, is further developed and applied to
the study of multiple zeta values.  In particular, we establish
evaluations for certain sums of cyclically generated multiple zeta
values.  The boundary case of our result reduces to a former
conjecture of Zagier.}


\newcommand\shuff{
\setlength{\unitlength}{.4pt}
\begin{picture}(40,20)
\put(10,2){\line(1,0){20}} \put(10,2){\line(0,1){10}}
\put(20,2){\line(0,1){10}} \put(30,2){\line(0,1){10}}
\end{picture}
}

\newcommand\shuffwithlimits{\setlength{\unitlength}{.6pt}
\begin{picture}(40,20)
\put(10,2){\line(1,0){20}} \put(10,2){\line(0,1){15}}
\put(20,2){\line(0,1){15}} \put(30,2){\line(0,1){15}}
\end{picture}
}

\newcommand{\iso}{\stackrel{\sim}{\to}}
\newcommand{\iu}{\int_0^1}
\newcommand{\W}{{\mathbf W}}
\newcommand{\X}{{\mathbf X}}
\newcommand{\Y}{{\mathbf Y}}
\newcommand{\Z}{{\mathbf Z}}
\newcommand{\R}{{\mathbf R}}
\newcommand{\Q}{{\mathbf Q}}
\newcommand{\C}{{\mathbf C}}
\newcommand{\Sh}{{\mathrm{Sh}}}
\newcommand{\Shuff}{{\mathrm{Shuff}}}
\newcommand{\ShuffProd}[2]{\mathop{\shuffwithlimits}\limits_{#1}^{#2}}
\newcommand{\z}{\zeta}

\newtheorem{Conj}{Conjecture}

\begin{article}


\section{Introduction}
\label{sect:motive}

We continue our study of nested sums of the form
\begin{equation}
   \z(s_1,s_2,\dots,s_k):=\sum_{n_1>n_2>\cdots>n_k>0}
   \;\prod_{j=1}^k n_j^{-s_j},
\label{MZVdef}
\end{equation}
commonly referred to as  multiple zeta
values~\cite{BBBLa,BBBLc,BowBrad,Hoff4,Hoff5,YOhno,Zag}.  Here and
throughout, $s_1,s_2,\dots,s_k$ are positive integers with $s_1>1$
to ensure convergence.

There exist many intriguing results and conjectures concerning
values of~(\ref{MZVdef}) at various arguments.  For example,
\begin{equation}
   \z(\{3,1\}^n) :=
   \z(\underbrace{3,1,3,1,\ldots,3,1}_{2n})
   = \frac{2\pi^{4n}}{(4n+2)!},
   \qquad 0\le n\in\Z,
\label{Z31}
\end{equation}
was conjectured by Zagier~\cite{Zag} and first proved by
Broadhurst et al~\cite{BBBLa} using analytic techniques.
Subsequently, a purely combinatorial proof was given~\cite{BBBLc}
based on the well-known shuffle property of iterated integrals,
and it is this latter approach which we develop more fully here.
For further and deeper results from the analytic viewpoint,
see~\cite{BowBrad}.

Our main result is a generalization of~(\ref{Z31}) in which twos
are inserted at various places in the argument string $\{3,1\}^n$.
Given a non-negative integer $n$, let $\vec
s=(m_0,m_1,\dots,m_{2n})$ be a vector of non-negative integers,
and consider the multiple zeta value obtained by inserting $m_j$
consecutive twos after the $j$th element of the string $\{3,1\}^n$
for each $j=0,1,2,\dots,2n$:
\begin{multline*}
   Z(\vec s)\\
   :=\z(\{2\}^{m_0},3,\{2\}^{m_1},1,\{2\}^{m_2},3,\{2\}^{m_3},1,
   \dots,3,\{2\}^{m_{2n-1}},1,\{2\}^{m_{2n}}).
\end{multline*}
For non-negative integers $k$ and $r$, let $C_r(k)$ denote the set
of $\binom{k+r-1}{r-1}$ ordered non-negative integer compositions
of $k$ having $r$ parts. For example,
$C_3(2)=\{(2,0,0),(0,2,0),(0,0,2),(0,1,1),(1,0,1),(1,1,0)\}$. Our
generalization of~(\ref{Z31}) states (see Corollary~\ref{cor:T} of
Section~\ref{sect:cycle}) that
\begin{equation}
   \sum_{\vec s\in C_{2n+1}(m-2n)} Z(\vec s) = \frac{2\pi^{2m}}{(2m+2)!}
   \binom{m+1}{2n+1},
\label{Zsum}
\end{equation}
for all non-negative integers $m$ and $n$ with $m\ge 2n$.
Equation~(\ref{Z31}) is the special case of~(\ref{Zsum}) in which
$m=2n$, since $Z(\{0\}^{2n+1})=\z(\{3,1\}^n)$.  If again $\vec
s=(m_0,m_1,\dots,m_{2n})$ and we put
\[
   {\mathscr{C}}(\vec s) :=Z(\vec s)+\sum_{j=1}^{2n}
   Z(m_j,m_{j+1},\dots,m_{2n},m_0,\dots,m_{j-1}),
\]
then (see Theorem~\ref{thm:T} of Section~\ref{sect:cycle})
\begin{equation}
   \sum_{\vec s\in C_{2n+1}(m-2n)}{\mathscr{C}}(\vec s) =
   Z(m)\times|C_{2n+1}(m-2n)| =
   \frac{\pi^{2m}}{(2m+1)!}\binom{m}{2n}
\label{Csum}
\end{equation}
is an equivalent formulation of~(\ref{Zsum}).  The cyclic
insertion conjecture~\cite{BBBLc} can be restated as the assertion
that ${\mathscr{C}}(\vec s)=Z(m)$ for all $\vec s\in
C_{2n+1}(m-2n)$ and integers $m\ge 2n\ge 0$.  Thus, our result
reduces the problem to that of establishing the invariance of
${\mathscr{C}}(\vec s)$ on $C_{2n+1}(m-2n)$.

The outline of the paper is as follows.
Section~\ref{sect:IterInts} provides the essential background for
our results.  The theory is formalized and further developed in
Section~\ref{sect:ShuffAlg}, in which we additionally give a
simple proof of Ree's formula for the inverse of a Lie
exponential. In Section~\ref{sect:ShuffBasis} we focus on the
combinatorics of two-letter words, as this is most directly
relevant to the study of multiple zeta values.  In the final
section, we establish the aforementioned results~(\ref{Zsum})
and~(\ref{Csum}).

\section{Iterated Integrals}\label{sect:IterInts}
As Kontsevich~\cite{Zag} observed,~(\ref{MZVdef}) admits an
iterated integral representation
\begin{equation}
   \z(s_1,s_2,\dots,s_k) = \iu \prod_{j=1}^k a^{s_j-1}b
\label{iterint}
\end{equation}
of depth $\sum_{j=1}^k s_j$.  Here, the notation
\begin{equation}
   \int_y^x \prod_{j=1}^n \alpha_j :=
   \Int_{x>t_1>t_2>\cdots>t_n>y} \;\prod_{j=1}^n f_j(t_j)\,dt_j,
   \qquad \alpha_j := f_j(t_j)\,dt_j
\label{IterIntNotn}
\end{equation}
of~\cite{BBBLa} is used with $a$ and $b$ denoting the differential
1-forms $dt/t$ and $dt/(1-t)$, respectively.  Thus, for example,
if $f_1\ne f_2$, we write $\alpha_1^2\alpha_2\alpha_1$ for the
integrand
$f_1(t_1)f_1(t_2)f_2(t_3)f_1(t_4)\,dt_1\,dt_2\,dt_3\,dt_4.$
Furthermore, we shall agree that any iterated integral of an empty
product of differential $1$-forms is equal to $1$.  This
convention is mainly a notational convenience; nevertheless we
shall find it useful for stating results about iterated integrals
more concisely and naturally than would be possible otherwise.
Thus~(\ref{IterIntNotn}) reduces to $1$ when $n=0$ regardless of
the values of $x$ and $y$.

Clearly the product of two iterated integrals of the
form~(\ref{IterIntNotn}) consists of a sum of iterated integrals
involving all possible interlacings of the variables. Thus if we
denote the set of all $\binom{n+m}{n}$ permutations $\sigma$ of
the indices $\{1,2,\dots,n+m\}$ satisfying
$\sigma^{-1}(j)<\sigma^{-1}(k)$ for all $1\le j<k\le n$ and
$n+1\le j<k\le n+m$ by $\Shuff(n,m)$, then we have the
self-evident formula
\[
   \bigg(\int_y^x \prod_{j=1}^n \alpha_j \bigg)\bigg(\int_y^x\,
   \prod_{j=n+1}^{n+m}\alpha_j\bigg)
   = \sum_{\sigma\in\Shuff(n,m)} \int_y^x\; \prod_{j=1}^{n+m}
   \alpha_{\sigma(j)},
\]
and so define the shuffle product $\shuff$ by
\begin{equation}
\label{shuff-rule}
   \bigg(\prod_{j=1}^n \alpha_j\bigg)\shuff
   \bigg(\prod_{j=n+1}^{n+m}\alpha_j\bigg)
    := \sum_{\sigma\in\Shuff(n,m)} \prod_{j=1}^{n+m}
    \alpha_{\sigma(j)}.
\end{equation}
Thus, the sum is over all non-commutative products (counting
multiplicity) of length $n+m$ in which the relative orders of the
factors in the products $\alpha_1\alpha_2\cdots \alpha_n$ and
$\alpha_{n+1}\alpha_{n+2}\cdots \alpha_{n+m}$ are preserved.  The
term ``shuffle'' is used because such permutations arise in riffle
shuffling a deck of $n+m$ cards cut into one pile of $n$ cards and
a second pile of $m$ cards.

The study of shuffles and iterated integrals was pioneered by
Chen~\cite{Chen54,Chen57} and subsequently formalized by
Ree~\cite{Ree}.  A fundamental formula noted by Chen expresses an
iterated integral of a product of two paths as a convolution of
iterated integrals over the two separate paths.  A second formula
also due to Chen shows what happens when the underlying
simplex~(\ref{IterIntNotn}) is re-oriented.   Chen's proof in both
cases is by induction on the number of differential $1$-forms.
Since we will make use of these results in the sequel, it is
convenient to restate them here in the current notation and give
direct proofs.

\begin{proposition}[{\cite[(1.6.2)]{Chen71}}]\label{prop:ReverseForms}
Let $\alpha_1,\alpha_2,\dots,\alpha_n$ be differential
$1$-\linebreak[4]forms and let $x,y\in\R$.  Then
\[
   \int_y^x \alpha_1\alpha_2\cdots\alpha_n
   = (-1)^n \int_x^y \alpha_n\alpha_{n-1}\cdots\alpha_1.
\]
\end{proposition}

\begin{proof}
Suppose $\alpha_j=f_j(t_j)\,dt_j$.  Observe that
\begin{eqnarray*}
   &&\int_y^x f_1(t_1)\int_y^{t_1}f_2(t_2)\cdots\int_y^{t_{n-1}}
   f_n(t_n)\,dt_n\,dt_{n-1}\cdots\,dt_1\\
   &&\qquad =\int_y^x f_n(t_n)\int_{t_n}^x
   f_{n-1}(t_{n-1})\cdots\int_{t_2}^x
   f_1(t_1)\,dt_1\,dt_2\cdots\,dt_n.
\end{eqnarray*}
Now switch the limits of integration at each level.
\end{proof}

\begin{proposition}[{\cite[Lemma 1.1]{Chen54}}]\label{prop:HC}
Let $\alpha_1,\alpha_2,\dots,\alpha_n$ be differential $1$-forms
and let $y\le z\le x$.  Then
\[
   \int_y^x \prod_{j=1}^n \alpha_j
   = \sum_{k=0}^n
     \bigg(\int_z^x \prod_{j=1}^k \alpha_j\bigg)
     \bigg(\int_y^z \prod_{j=k+1}^n \alpha_j\bigg).
\]
\end{proposition}

\begin{proof}
\begin{eqnarray*}
   &&\big\{(t_1,t_2,\dots,t_n)\in\R^n :
   x>t_1>t_2>\cdots>t_n>y\big\}\\
   &&\qquad = \bigcup_{k=0}^n \big\{(t_1,\dots,t_k)\in\R^k :
   x>t_1>\cdots>t_k>z\big\}\\
   &&\qquad\qquad\qquad\times
   \big\{(t_{k+1},\dots,t_n)\in\R^{n-k} : z>t_{k+1}>\cdots>t_n>y\big\}.
\end{eqnarray*}
\end{proof}

A related version of Proposition~\ref{prop:HC}, ``H\"older
Convolution,'' is exploited in~\cite{BBBLa} to indicate how rapid
computation of multiple zeta values and related slowly-convergent
multiple polylogarithmic sums is accomplished.  In
Section~\ref{sect:Lie}, Proposition~\ref{prop:HC} is used in
conjunction with Proposition~\ref{prop:ReverseForms} to give a
quick proof of Ree's formula~\cite{Ree} for the inverse of a Lie
exponential.

\section{The Shuffle Algebra}
\label{sect:ShuffAlg}

We have seen how shuffles arise in the study of iterated integral
representations for multiple zeta values. Following~\cite{MinPet1}
(cf.\ also~\cite{BBBLc,Ree}) let $A$ be a finite set and let $A^*$
denote the free monoid generated by $A$. We regard $A$ as an
alphabet, and the elements of $A^*$ as words formed by
concatenating any finite number of letters from this alphabet.  By
linearly extending the concatenation product to the set $\Q\langle
A\rangle$ of rational linear combinations of elements of $A^*$, we
obtain a non-commutative polynomial ring with indeterminates the
elements of $A$ and with multiplicative identity $1$ denoting the
empty word.

The shuffle product is alternatively defined first on words by the
recursion
\begin{equation}
\begin{cases}
   \forall w\in A^*, \quad & 1\shuff w = w\shuff 1 = w,\\
   \forall a,b\in A, \quad\forall u,v\in A^*, \quad & au\shuff bv
   =a(u\shuff bv)+b(au\shuff v),
\end{cases}
\label{LeftShuffDef}
\end{equation}
and then extended linearly to $\Q\langle A\rangle$.  One checks
that the shuffle product so defined is associative and
commutative, and thus $\Q\langle A\rangle$ equipped with the
shuffle product becomes a commutative $\Q$-algebra, denoted
$\Sh_{\Q}[A]$.  Radford~\cite{Rad} has shown that $\Sh_{\Q}[A]$ is
isomorphic to the polynomial algebra $\Q[L]$ obtained by adjoining
to $\Q$ the transcendence basis $L$ of Lyndon words.

The recursive definition~(\ref{LeftShuffDef}) has its analytical
motivation in the formula for integration by parts---equivalently,
the product rule for differentiation.  Thus, if we put
$a=f(t)\,dt$, $b=g(t)\,dt$ and
\[
   F(x) := \int_y^x (au\shuff bv) =
   \bigg(\int_y^x f(t)\int_y^t u\,dt\bigg)\bigg(\int_y^x g(t)\int_y^t
   v\,dt\bigg),
\]
then writing $F(x)=\int_y^x F'(s)\,ds$ and applying the product
rule for differentiation yields
\begin{eqnarray*}
   F(x) &=& \int_y^x \bigg(f(s)\int_y^s u\bigg)\bigg(\int_y^s
   g(t)\int_y^t v\,dt\bigg)\,ds\\
   && \qquad\qquad + \int_y^x g(s)\bigg(\int_y^s
   f(t)\int_y^t u\,dt\bigg)\int_y^s v\,ds\\
   &=& \int_y^x \left[ a(u\shuff bv) + b(au\shuff v)\right].
\end{eqnarray*}
Alternatively, by viewing $F$ as a function of $y$, we see that
the recursion could equally well have been stated as
\begin{equation}
\begin{cases}
   \forall w\in A^*, \quad & 1\shuff w = w\shuff 1 = w,\\
   \forall a,b\in A, \quad\forall u,v\in A^*, \quad & ua\shuff vb
   =(u\shuff vb)a+(ua\shuff v)b.
\end{cases}
\label{RightShuffDef}
\end{equation}
Of course, both definitions are equivalent to~(\ref{shuff-rule}).

\subsection{$\Q$-Algebra Homomorphisms on Shuffle
Algebras}\label{subsect:hom}

The following relatively straightforward results concerning
$\Q$-algebra homomorphisms on shuffle algebras will facilitate our
discussion of the Lie exponential in Section~\ref{sect:Lie} and of
relationships between certain identities for multiple zeta values
and Euler sums~\cite{BBB,BBBLa,BowBrad}.  To reduce the
possibility of any confusion in what follows, we make the
following definition explicit.

\begin{definition}\label{def:AntiHom}
Let $R$ and $S$ be rings with identity, and let $A$ and $B$ be
alphabets.  A ring anti-homomorphism $\psi : R\langle A\rangle\to
S\langle B\rangle$ is an additive, $R$-linear, identity-preserving
map that satisfies $\psi(u)\psi(v)=\psi(vu)$ for all $u,v\in A^*$
(and hence for all $u,v\in R\langle A\rangle$).
\end{definition}

\begin{proposition}\label{prop:AntiHom} Let $A$ and $B$ be alphabets.
A ring anti-homomor\-phism $\psi:\Q\langle A\rangle\to \Q\langle
B\rangle$ that satisfies $\psi(A)\subseteq B$ induces a
$\Q$-algebra homomorphism of shuffle algebras $\psi:\Sh_{\Q}[A]
\to\Sh_{\Q}[B]$ in the natural way.
\end{proposition}

\begin{proof}  It suffices to show that
$\psi(u\shuff v)=\psi(u)\shuff\psi(v)$ for all $u,v\in A^*$. The
proof is by induction, and will require both recursive definitions
of the shuffle product.  Let $u,v\in A^*$ be words. For the base
case, note that $\psi(1\shuff u)=\psi(u)=1\shuff\psi(u)$ and
likewise with the empty word on the right.  For the inductive
step, let $a,b\in A$ be letters and assume that $\psi(u\shuff
bv)=\psi(u)\shuff \psi(bv)$ and $\psi(au\shuff
v)=\psi(au)\shuff\psi(v)$ both hold.  Then as $\psi$ is an
anti-homomorphism of rings,
\begin{eqnarray*}
   \psi(au\shuff bv) &=& \psi(a(u\shuff bv)+b(au\shuff v))\\
   &=&\psi(a(u\shuff bv))+\psi(b(au\shuff v))\\
   &=&\psi(u\shuff bv)\psi(a)+\psi(au\shuff v)\psi(b)\\
   &=&[\psi(u)\shuff\psi(bv)]\psi(a)+[\psi(au)\shuff\psi(v)]\psi(b)\\
   &=&[\psi(u)\shuff\psi(v)\psi(b)]\psi(a)+[\psi(u)\psi(a)\shuff\psi(v)]
      \psi(b)\\
   &=&\psi(u)\psi(a)\shuff \psi(v)\psi(b)\\
   &=&\psi(au)\shuff\psi(bv).
\end{eqnarray*}
\end{proof}

Of course, there is an analogous result for ring homomorphisms.

\begin{proposition}\label{prop:Hom} Let $A$ and $B$ be alphabets.  A ring
homomorphism $\phi:\Q\langle A\rangle\to\Q\langle B\rangle$ that
satisfies $\phi(A)\subseteq B$ induces a $\Q$-algebra homomorphism
of shuffle algebras $\phi:\Sh_{\Q}[A]\to \Sh_{\Q}[B]$ in the
natural way.
\end{proposition}

\begin{proof}  The proof is similar to the proof of
Proposition~\ref{prop:AntiHom}, and in fact is simpler in that it
requires only one of the two recursive definitions of the shuffle
product.  Alternatively, one can put $u=a_1 a_2\cdots a_n$,
$v=a_{n+1} a_{n+2}\cdots a_{n+m}$ and verify the equation
$\phi(u\shuff v)=\phi(u)\shuff\phi(v)$ using~(\ref{shuff-rule})
and the hypothesis that $\phi$ is a ring homomorphism on
$\Q\langle A\rangle$. \end{proof}

\begin{demo} {Example 1} Let $A$ be an alphabet and
let $R:\Q\langle A\rangle \to\Q\langle A\rangle$ be the canonical
ring anti-automorphism induced by the assignments $R(a)=a$ for all
$a\in A$.  Then $R(\prod_{j=1}^n a_j)=\prod_{j=1}^n a_{n-j+1}$ for
all $a_1,\dots,a_n\in A$, so that $R$ is a string-reversing
involution which induces a shuffle algebra automorphism of
$\Sh_{\Q}[A]$.  We shall reserve the notation $R$ for this
automorphism throughout.
\end{demo}

\begin{demo} {Example 2} Let $A=\{a,b\}$ and let $S:\Q\langle
A\rangle \to\Q\langle A\rangle$ be the ring automorphism induced
by the assignments $S(a)=b$, $S(b)=a.$ Then the composition $\psi
:= S\circ R$ is a letter-switching, string-reversing involution
which induces a shuffle algebra automorphism of $\Sh_{\Q}[A]$. In
the case $a=dt/t$, $b=dt/(1-t)$, this is the so-called Kontsevich
duality~\cite{Zag,BBB,BBBLa,YOhno} for iterated integrals obtained
by making the change of variable $t\mapsto 1-t$ at each level of
integration.  Words which are invariant under $\psi$ are referred
to as \textit{self-dual}.  It is easy to see that a self-dual word
must be of even length, and the number of self-dual words of
length $2k$ is $2^k$.
\end{demo}

\begin{demo} {Example 3} Let $A=\{a,b\}, B=\{b,c\}$ and let
$\psi:\Q\langle A\rangle \to \Q\langle B\rangle$ be the
letter-shifting, string-reversing ring anti-homomorphism induced
by the assignments $\psi(a)=b$ and $\psi(b)= c$. Then $\psi$
induces a shuffle algebra isomorphism $\psi:\Sh_{\Q}[A]\iso
\Sh_{\Q}[B]$. With the choice of differential $1$-forms $a=dt/t$,
$b=dt/(1-t)$, $c=-dt/(1+t)$, $\psi$ maps shuffle identities for
multiple zeta values to equivalent identities for alternating unit
Euler sums. We refer the reader to~\cite{BBB,BBBLa,BowBrad} for
details concerning alternating Euler sums; for our purposes here
it suffices to assert that they are important instances---as are
multiple zeta values---of multiple
polylogarithms~\cite{BBBLa,Gonch}.
\end{demo}

\subsection{A Lie Exponential}\label{sect:Lie}
Let $A$ be an alphabet, and let $X=\{X_a : a\in A\}$ be a set of
${\mathrm{card}}(A)$ distinct non-commuting indeterminates.  Every
element in $\Q\langle X\rangle$ can be written as a sum
$F=F_0+F_1+\cdots$ where $F_n$ is a homogeneous form of degree
$n$. Those elements $F$ for which $F_n$ belongs to the Lie algebra
generated by $X$ for each $n>0$ and for which $F_0=0$ are referred
to as {\em Lie elements}.

Let $\X:\Q\langle A\rangle \to \Q\langle X\rangle$ be the
canonical ring isomorphism induced by the assignments $\X(a)=X_a$
for all $a\in A$.  If $Y=\{Y_a: a\in A\}$ is another set of
non-commuting indeterminates, we similarly define $\Y:\Q\langle
A\rangle \to \Q\langle Y\rangle$ to be the canonical ring
isomorphism induced by the assignments $\Y(a)=Y_a$ for all $a\in
A$. Let us suppose $X=\X(A)$ and $Y=\Y(A)$ are disjoint and their
elements commute with each other, so that for all $a,b\in A$ we
have $X_aY_b=Y_bX_a$. If we define addition and multiplication in
$\Q[\X,\Y]$ by $(\X+\Y)(a)=X_a+Y_a$ and $(\X\Y)(a)=X_aY_a$ for all
$a\in A$, then $\Q[\X,\Y]$ becomes a commutative $\Q$-algebra of
ring isomorphisms $\Z$. For example, if $\Z=\X+\Y$ and
$w=a_1a_2\cdots a_n$ where $a_1,a_2\dots,a_n\in A$, then
\[
   \Z(w) = (\X+\Y)(a_1a_2\cdots a_n)
           = \prod_{j=1}^n (\X+\Y)(a_j)
           = \prod_{j=1}^n \left(X_{a_j}+Y_{a_j}\right).
\]
Let $G:\Q[\X,\Y] \to (\Sh_{\Q}[A])\langle\langle
X,Y\rangle\rangle$ be defined by
\begin{equation}
   G(\Z) := \sum_{w\in A^*} w \Z(w).
\label{Gdef}
\end{equation}
Evidently,
\begin{equation}
   G(\X) = 1+\sum_{n=1}^\infty\bigg(\sum_{a\in A} a X_a\bigg)^n =
   \frac{1}{1-\sum_{a\in A}a X_a}.
\label{Gogf}
\end{equation}
More importantly, $G$ is a homomorphism from the underlying
$\Q$-vector space to the underlying multiplicative monoid
$((\Sh_{\Q}[A])\langle\langle X,Y\rangle \rangle,\shuff).$

\begin{theorem}\label{thm:PlusToShuff}  The map
$G:\Q[\X,\Y]\to (\Sh_{\Q}[A])\langle\langle X,Y\rangle \rangle$
defined by (\ref{Gdef}) has the property that
\[
   G(\X+\Y) = G(\X)\shuff G(\Y).
\]
\end{theorem}

\begin{proof}
On the one hand, we have
\[
   G(\X+\Y) = \sum_{w\in A^*} w (\X+\Y)(w),
\]
whereas on the other hand,
\[
   G(\X)\shuff G(\Y) = \sum_{u\in A^*} u\X(u) \shuff\sum_{v\in A^*}
   v \Y(v) = \sum_{u,v\in A^*}(u\shuff v) \X(u)\Y(v).
\]
Therefore, we need to show that
\[
   \sum_{u,v\in A^*} (u\shuff v) \X(u)\Y(v) = \sum_{w\in A^*}w(\X+\Y)(w).
\]
But,
\begin{eqnarray*}
   &&\sum_{u,v\in A^*} (u\shuff v) \X(u)\Y(v)\\
   &=& \sum_{n\ge 0} \; \sum_{a_1,\dots,a_n\in A}\;
   \sum_{k=0}^n \bigg(\prod_{j=1}^k a_j \shuff \prod_{j=k+1}^n
   a_j\bigg) \prod_{j=1}^k X_{a_j}\prod_{j=k+1}^n
   Y_{a_j}\\
   &=& \sum_{n\ge 0}\;\sum_{a_1,\dots,a_n\in A}\;
   \sum_{k=0}^n \;\sum_{\sigma\in \Shuff(k,n-k)}\;
   \prod_{r=1}^n a_{\sigma(r)} \prod_{j=1}^k X_{a_j}
   \prod_{j=k+1}^n Y_{a_j},
\end{eqnarray*}
using the non-recursive definition~(\ref{shuff-rule}) of the
shuffle product. For each $\sigma\in \Shuff(k,n-k)$, if
$a_1,\dots,a_n$ run through the elements of $A$, then so do
$a_{\sigma(1)},\dots,a_{\sigma(n)}$.  Hence putting
$b_j=a_{\sigma(j)}$, we have that
\begin{eqnarray*}
   &&\sum_{u,v\in A^*} (u\shuff v) \X(u)\Y(v)\\
   &=&\sum_{n\ge 0}\;\sum_{b_1,\dots,b_n\in A}\;
   \bigg(\prod_{r=1}^n b_r\bigg)\sum_{k=0}^n\;
   \sum_{\sigma\in\Shuff(k,n-k)}\; \prod_{j=1}^k X_{b_{\sigma^{-1}(j)}}
   \prod_{j=k+1}^n Y_{b_{\sigma^{-1}(j)}}\\
   &=& \sum_{n\ge 0}\;\sum_{b_1,\dots,b_n\in A}\;
   \bigg(\prod_{r=1}^n b_r\bigg)
   \prod_{j=1}^n (X_{b_j}+Y_{b_j})\\
   &=& \sum_{w\in A^*} w(\X+\Y)(w).
\end{eqnarray*}
In the penultimate step, we have summed over all $\binom{n}{k}$
shuffles of the indeterminates $X$ with the indeterminates $Y,$
yielding all $2^n$ possible choices obtained by selecting an $X$
or a $Y$ from each factor in the product
$(X_{b_1}+Y_{b_1})\cdots(X_{b_n}+Y_{b_n})$. \end{proof}

\begin{demo} {Remarks} Theorem~\ref{thm:PlusToShuff} suggests that
the map $G$ defined by~(\ref{Gdef}) can be viewed as a
non-commutative analog of the exponential function.  The analogy
is clearer if we rewrite~(\ref{Gogf}) in the form
\[
   G(\X) = 1+ \sum_{n=1}^\infty \frac{1}{n!}\bigg(\sum_{a\in A} a
   X_a\bigg)^{\!\shuff n}.
\]
Just as the functional equation for the exponential function is
equivalent to the binomial theorem, Theorem~\ref{thm:PlusToShuff}
is equivalent to the following shuffle analog of the binomial
theorem:
\begin{proposition}[Binomial Theorem in $\Q\langle X\rangle \langle
Y\rangle$] Let $X=\{X_1,\linebreak[0] X_2,\dots,X_n\}$ and $Y=
\{Y_1,Y_2,\dots,Y_n\}$ be disjoint sets of non-commuting
indeterminates such that $X_jY_k=Y_kX_j$ for all $1\le j,k\le n$.
Then
\[
   \prod_{j=1}^n (X_j+Y_j) = \sum_{k=0}^n\;
   \sum_{\sigma\in\Shuff(k,n-k)}\;\prod_{j=1}^k X_{\sigma^{-1}(j)}
   \prod_{j=k+1}^n Y_{\sigma^{-1}(j)}.
\]
\end{proposition}
\end{demo}

Chen~\cite{Chen54,Chen57} considered what is in our notation the
iterated integral of~(\ref{Gdef}), namely
\begin{equation}
   G_y^x:=   \sum_{w\in A^*} \int_y^x w \X(w)
\label{ChenSeries}
\end{equation}
in which the alphabet $A$ is viewed as a set of differential
$1$-forms.  He
proved~\cite[Theorem~6.1]{Chen54},~\cite[Theorem~2.1]{Chen57} the
non-commutative generating function formulation
\[
   G_y^x = G_z^x G_y^z,\qquad y\le z\le x
\]
of Proposition~\ref{prop:HC} and also
proved~\cite[Theorem~4.2]{Chen57} that if the $1$-forms are
piecewise continuously differentiable, then $\log G_y^x$ is a Lie
element, or equivalently, that $G_y^x$ is a Lie exponential.
However, Ree~\cite{Ree} showed that a formal power series
\[
   \log\bigg(1+\sum_{n>0}\; \sum_{1\le j_1,\dots,j_n\le m}
   c(j_1,\dots,j_n) X_{j_1}\cdots X_{j_n}\bigg)
\]
in non-commuting indeterminates $X_j$ is a Lie element if and only
if the coefficients satisfy the shuffle relations
\[
   c(j_1,\dots,j_n)c(j_{n+1},\dots,j_{n+k})
   =\sum_{\sigma\in\Shuff(n,k)}
   c(j_{\sigma(1)},\dots,j_{\sigma(n+k)}),
\]
for all non-negative integers $n$ and $k$. Using integration by
parts, Ree~\cite{Ree} showed that Chen's coefficients do indeed
satisfy these relations, and that more generally, $G(\X)$ as
defined by~(\ref{Gdef}) is a Lie exponential, a fact that can also
be deduced from Theorem~\ref{thm:PlusToShuff} and a result
of~Friedrichs~\cite{Fried,Lyn,Mag}.

Ree also proved a formula~\cite[Theorem~2.6]{Ree} for the inverse
of~(\ref{Gdef}), using certain derivations and Lie bracket
operations. It may be of interest to give a more direct proof,
using only the shuffle operation. The result is restated below in
our notation.
\begin{theorem}[{\cite[Theorem~2.6]{Ree}}]\label{thm:Ree}
Let $A$ be an alphabet, let $X=\{X_a : a\in A\}$ be a set of
non-commuting indeterminates and let $\X:\Q\langle
A\rangle\to\Q\langle X\rangle$ be the canonical ring isomorphism
induced by the assignments $\X(a)=X_a$ for all $a\in A$.   Let
$G(\X)$ be as in~(\ref{Gogf}), let $R$ be as in Example~1, and put
\[
   H(\X) := \sum_{w\in A^*} (-1)^{|w|}R(w) \X(w),
\]
where $|w|$ denotes the length of the word $w$.  Then $G(\X)\shuff
H(\X)=1$.
\end{theorem}
It is convenient to state the essential ingredient in our proof of
Theorem~\ref{thm:Ree} as an independent result.
\begin{lemma}\label{lem:Vacuum} Let $A$ be an
alphabet and let $R$ be as in Example~1.  For all $w\in A^*$, we
have
\begin{equation}
   \sum_{\substack{u,v\in A^*\\uv=w}}(-1)^{|u|}R(u)\shuff v =
   \delta_{|w|,0}.
\label{VacuumFront-Front}
\end{equation}
\end{lemma}

\begin{demo} {Remarks} We have used the Kronecker delta
\[
    \delta_{n,k}:= \begin{cases}
                    1 &\mbox{if $n=k$},\\
                    0 &\mbox{otherwise.}
                    \end{cases}
\]
Since $R$ is a $\Q$-algebra automorphism of $\Sh_{\Q}[A]$,
applying $R$ to both sides of~(\ref{VacuumFront-Front}) yields the
related identity
\[
   \sum_{\substack{u,v\in A^*\\ uv=w}}(-1)^{|u|}\,u\shuff R(v)=
   \delta_{|w|,0}, \qquad w\in A^*.
\]
\end{demo}

\begin{demo} {Proof of Lemma~\ref{lem:Vacuum}}
First note that if we view the elements of $A$ as differential
$1$-forms and integrate the left hand side
of~(\ref{VacuumFront-Front}) from $y$ to $x$, then we obtain
\[
   \sum_{\substack{u,v\in A^*\\uv=w}} (-1)^{|u|} \int_y^x
   R(u)\int_y^x v = \sum_{\substack{u,v\in A^*\\uv=w}}
   \int_x^y u\int_y^x v = \int_y^y w = \delta_{|w|,0}
\]
by Propositions~\ref{prop:ReverseForms} and~\ref{prop:HC}. For an
integral-free proof, we proceed as follows.
Clearly~(\ref{VacuumFront-Front}) holds when $|w|=0$, so assume
$w=\prod_{j=1}^n a_j$ where $a_1,\dots,a_n\in A$ and $n$ is a
positive integer.  Let ${\mathfrak{S}}_{n}$ denote the group of
permutations of the set of indices $\{1,2,\dots,n\}$, and let the
additive weight-function $W:2^{{\mathfrak{S}}_n}\to A^*$ map
subsets of ${\mathfrak{S}}_n$ to words as follows:
\[
   W(S) := \sum_{\sigma\in S} \prod_{j=1}^n a_{\sigma(j)},
   \qquad S\subseteq {\mathfrak{S}}_n.
\]
For $k=0,1,\dots,n$ let
\begin{eqnarray*}
   c_k &:=&
   W(\{\sigma\in{\mathfrak{S}}_n:\sigma^{-1}(i)<\sigma^{-1}(j)
      \;\mbox{for}\; k\ge i>j\ge 1\;\mbox{and}\\
      &&\qquad k+1\le i<j\le n\}),\\
   b_k &:=&
   W(\{\sigma\in{\mathfrak{S}}_n:\sigma^{-1}(i)<\sigma^{-1}(j)
      \;\mbox{for}\; k\ge i>j\ge 1\;\mbox{and}\\
      &&\qquad k\le i<j\le n\}).
\end{eqnarray*}
Then $c_0=b_1$, $c_n=b_n$ and $c_k=b_k+b_{k+1}$ for $1\le k\le
n-1$.  Thus,
\begin{eqnarray*}
   \sum_{\substack{u,v\in A^*\\uv=w}} (-1)^{|u|} R(u)\shuff v
   &=& \sum_{k=0}^n (-1)^k \prod_{j=1}^k a_{k-j+1}\shuff
   \prod_{j=k+1}^n a_j\\
   &=& \sum_{k=0}^n (-1)^k c_k\\
   &=& b_1+(-1)^n b_n + \sum_{k=1}^{n-1} (-1)^k (b_k+b_{k+1})\\
   &=& b_1+(-1)^n b_n +\sum_{k=1}^{n-1}(-1)^k b_k
      - \sum_{k=2}^n (-1)^k b_k\\
   &=& 0,
\end{eqnarray*}
since the sums telescope. $\qed$
\end{demo}

{\it Remark.} One can also give an integral-free proof of
Lemma~\ref{lem:Vacuum} by induction using the recursive
definition~(\ref{RightShuffDef}) of the shuffle product.

\begin{demo} {Proof of Theorem~\ref{thm:Ree}}
By Lemma~\ref{lem:Vacuum}, we have
\begin{eqnarray*}
   &&\sum_{u\in A^{*}} (-1)^{|u|}R(u)\X(u) \shuff \sum_{v\in
   A^*}v\X(v)\\
   &=&\sum_{w\in A^*} \X(w)\sum_{\substack{u,v\in
   A^*\\uv=w}}(-1)^{|u|}R(u)\shuff v\\
   &=&\sum_{w\in A^*} \X(w)\,\delta_{|w|,0}\\
   &=& 1.
\end{eqnarray*}
Since $(\Sh_{\Q}[A])\langle\langle X\rangle\rangle$ is commutative
with respect to the shuffle product, the result follows. $\qed$
\end{demo}

\section{Combinatorics of Shuffle Products}
\label{sect:ShuffBasis}

The combinatorial proof~\cite{BBBLc} of Zagier's
conjecture~(\ref{Z31}) hinged on expressing the sum of the words
comprising the shuffle product of $(ab)^p$ with $(ab)^q$ as a
linear combination of basis subsums $T_{p+q,n}$.  To gain a deeper
understanding of the combinatorics of shuffles on two letters, it
is necessary to introduce additional basis subsums. We do so here,
and thereby find analogous expansion theorems.  We conclude the
section by providing generating function formulations for these
results.  The generating function formulation plays a key role in
the proof of our main result~(\ref{Csum}), Theorem~\ref{thm:T} of
Section~\ref{sect:cycle}. The precise definitions of the basis
subsums follow.
\begin{definition}(\cite{BBBLc})\label{def:T}
For integers $m\ge n\ge 0$ let $S_{m,n}$ denote the set of words
occurring in the shuffle product $(ab)^{n}\shuff (ab)^{m-n}$ in
which the subword $a^2$ appears exactly $n$ times, and let
$T_{m,n}$ be the sum of the $\binom{m}{2n}$ distinct words in
$S_{m,n}.$  For all other integer pairs $(m,n)$ it is convenient
to define $T_{m,n}:=0$.
\end{definition}

\begin{definition}\label{def:U} For integers $m\ge n+1\ge 2$, let
$U_{m,n}$ be the sum of the elements of the set of words arising
in the shuffle product of $b(ab)^{n-1}$ with $b(ab)^{m-n-1}$ in
which the subword $b^2$ occurs exactly $n$ times.  For all other
integer pairs $(m,n)$ define $U_{m,n}:=0.$
\end{definition}
In terms of the basis subsums, we have the following
decompositions:
\begin{proposition}[{\cite[Prop.~1]{BBBLc}}]
\label{prop:T-basis} For all non-negative integers $p$ and $q$,
\begin{equation}
   (ab)^p \shuff (ab)^q = \sum_{n=0}^{\min(p,q)} 4^n
   \binom{p+q-2n}{p-n}T_{p+q,n}.
\label{T-Basis}
\end{equation}
\end{proposition}
The corresponding result for our basis (Definition~\ref{def:U}) is
\begin{proposition}\label{prop:U-basis} For all positive integers $p$ and
$q$,
\begin{equation}
   b(ab)^{p-1} \shuff b(ab)^{q-1} = \frac12\sum_{n=1}^{\min(p,q)} 4^n
   \binom{p+q-2n}{p-n}U_{p+q,n}.
\label{U-Basis}
\end{equation}
\end{proposition}

\begin{demo} {Proof of Proposition~\ref{prop:U-basis}}  See the
proof of Proposition~\ref{prop:T-basis} given in~\cite{BBBLc}. The
only difference here is that $a^2$ occurs one less time per word
than $b^2$ and so the multiplicity of each word must be divided by
$2$. The index of summation now starts at $1$ because there must
be at least one occurrence of $b^2$ in each term of the expansion.
$\qed$
\end{demo}

\begin{corollary}
\label{Cor:comp} For integers $p\ge 1$ and  $q\ge 0$,
\begin{equation}
\label{combo}
\begin{split}
b(ab)^{p-1} \shuff (ab)^{q} &= \sum_{n=0}^{\min(p-1,q)} 4^n
   \binom{p+q-2n-1}{p-n-1}b\,T_{p+q-1,n}\\
&+ \frac12 \sum_{n=1}^{\min(p,q)} 4^n
   \binom{p+q-2n}{p-n}a\,U_{p+q,n}.
\end{split}
\end{equation}
\end{corollary}

\begin{proof}
>From (\ref{LeftShuffDef}) it is immediate that $$b(ab)^{p-1}
\shuff (ab)^{q}= b[(ab)^{p-1} \shuff (ab)^{q}] +a[b(ab)^{p-1}
\shuff b(ab)^{q-1}].$$ Now apply (\ref{T-Basis}) and Proposition
\ref{prop:U-basis}.\end{proof}

\begin{proposition}\label{prop:T-Umbral}
Let $x_0,x_1,\dots$ and $y_0,y_1,\dots$ be sequences of not
necessarily commuting indeterminates, and let $m$ be a
non-negative (respectively, positive) integer.
We have the shuffle convolution formulae
\begin{multline}
\label{T-Umbral}
   \sum_{k=0}^m
   x_k\, y_{m-k}
   \left[(ab)^k \shuff (ab)^{m-k}\right]\\
   = \sum_{n=0}^{\lfloor m/2\rfloor} 4^n\sum_{j=0}^{m-2n}
    \binom{m-2n}{j} x_{n+j}\, y_{m-n-j}\,T_{m,n},
\end{multline}
and
\begin{multline}
\label{U-Umbral}
   \sum_{k=1}^{m-1}
   x_k\, y_{m-k}
   \left[b(ab)^{k-1} \shuff b(ab)^{m-k-1}\right]\\
   = \frac12\sum_{n=1}^{\lfloor m/2\rfloor} 4^n\sum_{j=0}^{m-2n}
    \binom{m-2n}{j} x_{n+j}\, y_{m-n-j}\,U_{m,n},
\end{multline}
respectively.
\end{proposition}

\begin{proof}
Starting with the left hand side of~(\ref{T-Umbral}) and
applying~(\ref{T-Basis}), we find that
\begin{eqnarray*}
   && \sum_{k=0}^m x_k\, y_{m-k} \left[(ab)^k \shuff
   (ab)^{m-k}\right]\\
    &=& \sum_{k=0}^m x_k\, y_{m-k} \sum_{n=0}^{\min(k,m-k)}
    4^n \binom{m-2n}{k-n}T_{m,n}\\
    &=& \sum_{n=0}^{\lfloor m/2\rfloor} 4^n\,
    \sum_{k=n}^{m-n} x_k\, y_{m-k} \binom{m-2n}{k-n}\,T_{m,n}\\
    &=&  \sum_{n=0}^{\lfloor m/2\rfloor} 4^n\,
    \sum_{j=0}^{m-2n} \binom{m-2n}{j} x_{n+j}\, y_{m-n-j}\,T_{m,n},
\end{eqnarray*}
which proves~(\ref{T-Umbral}).  The proof of~(\ref{U-Umbral})
proceeds analogously from~(\ref{U-Basis}). \end{proof}

As the proof shows, the products taken in~(\ref{T-Umbral})
and~(\ref{U-Umbral}) can be quite general; between the not
necessarily commutative indeterminates and the polynomials in
$a,b$ the products need only be bilinear for the formul{\ae} to
hold. Thus, there are many possible special cases that can be
examined. Here we will consider only one major application.  If we
confine ourselves to commuting geometric sequences, we obtain
\begin{theorem}\label{thm:T-Binom}
Let $x$ and $y$ be commuting indeterminates.  In the commutative
polynomial ring $(\Sh_{\Q}[a,b])[x,y]$ we have the shuffle
convolution formulae
\begin{equation}
\label{T-Binom}
   \sum_{k=0}^m x^k y^{m-k} \left[(ab)^k \shuff (ab)^{m-k}\right]
   = \sum_{n=0}^{\lfloor m/2\rfloor} (4xy)^n (x+y)^{m-2n}\,
   T_{m,n}
\end{equation}
for all non-negative integers $m$, and
\begin{equation}
\label{U-Binom}
   \sum_{k=1}^{m-1} x^k y^{m-k} \left[b(ab)^{k-1} \shuff b(ab)^{m-k-1}\right]
   = \frac 12\sum_{n=1}^{\lfloor m/2\rfloor} (4xy)^n (x+y)^{m-2n}\, U_{m,n}
\end{equation}
for all integers $m\ge 2$.
\end{theorem}

\begin{proof}
In Proposition~\ref{prop:T-Umbral}, put $x_k=x^k$ and $y_k=y^k$
for each $k\ge 0$ and apply the binomial theorem. \end{proof}

\section{Cyclic Sums in $\Sh_{\Q}[a,b]$}
\label{sect:cycle}
\begin{sloppypar}
In this final section, we establish the results~(\ref{Zsum})
and~(\ref{Csum}) stated in the introduction.   Let $S_{m,n}$ be as
in Definition~\ref{def:T}. Each word in $S_{m,n}$ has a unique
representation
\begin{equation}
   (ab)^{m_0}\prod_{k=1}^n (a^2 b)(ab)^{m_{2k-1}}b(ab)^{m_{2k}},
\label{CompositionCorrespondence}
\end{equation}
in which $m_0,m_1,\dots,m_{2n}$ are non-negative integers with sum
$m_0+m_1+\cdots+m_{2n}=m-2n$.  Conversely, every ordered
$(2n+1)$-tuple $(m_0,m_1,\dots,m_{2n})$ of non-negative integers
with sum $m-2n$ gives rise to a unique word in $S_{m,n}$
via~(\ref{CompositionCorrespondence}).  Thus, a bijective
correspondence $\varphi$ is established between the set $S_{m,n}$
and the set $C_{2n+1}(m-2n)$ of ordered non-negative integer
compositions of $m-2n$ with $2n+1$ parts.  In view of the
relationship~(\ref{iterint}) expressing multiple zeta values as
iterated integrals, it therefore makes sense to define
\[
   Z(\vec s) := \iu \varphi(\vec s),\qquad \vec s\in C_{2n+1}(m-2n),
   \qquad a:= dt/t,\quad b:= dt/(1-t).
\]
Thus, if $\vec s=(m_0,m_1,\dots,m_{2n})$, then
\begin{eqnarray*}
   &&Z(\vec s)=
   \iu (ab)^{m_0}\prod_{k=1}^n (a^2b)(ab)^{m_{2k-1}}b(ab)^{m_{2k}}\\
   &=&\z(\{2\}^{m_0},3,\{2\}^{m_1},1,\{2\}^{m_2},3,\{2\}^{m_3},1,
   \dots,3,\{2\}^{m_{2n-1}},1,\{2\}^{m_{2n}}),
\end{eqnarray*}
in which the argument string consisting of $m_j$ consecutive twos
is inserted after the $j$th element of the string $\{3,1\}^n$ for
each $j=0,1,2,\dots,2n$.
\end{sloppypar}

\begin{sloppypar}
From~\cite{BBB} we recall the evaluation
\begin{equation}
   Z(m) = \z(\{2\}^m) = \frac{\pi^{2m}}{(2m+1)!},\qquad 0\le m\in\Z.
\label{z2repeated}
\end{equation}
Let ${\mathfrak{S}}_{2n+1}$ denote the group of permutations on
the set of indices $\{0,1,2,\dots,2n\}$. For
$\sigma\in{\mathfrak{S}}_{2n+1}$ we define a group action on
$C_{2n+1}(m-2n)$ by $\sigma\vec s =
(m_{\sigma^{-1}(0)},m_{\sigma^{-1}(1)},\dots,m_{\sigma^{-1}(2n)})$,
where $\vec s=(m_0,m_1,\dots,m_{2n})$.  Let
\begin{equation}
   {\mathscr{C}}(\vec s) := \sum_{j=0}^{2n} Z(\sigma^j\vec s),
   \qquad \sigma=(0\, 1\, 2 \dots 2n)
\label{CyclicSumDef}
\end{equation}
denote the sum of the $2n+1$ $Z$-values in which the arguments are
permuted cyclically.  By construction, ${\mathscr{C}}$ is
invariant under any cyclic permutation of its argument string. The
cyclic insertion conjecture~\cite[Conjecture 1]{BBBLc} asserts
that in fact, ${\mathscr{C}}$ depends only on the number and sum
of its arguments. More specifically, it is conjectured that
\begin{Conj}\label{conj:CyclicInsertion}
For any non-negative integers $m_0,m_1,\dots,m_{2n}$, we have
\[
   {\mathscr{C}}(m_0,m_1,\dots,m_{2n}) = Z(m)
   = \frac{\pi^{2m}}{(2m+1)!},
\]
where $m:=2n+\sum_{j=0}^{2n}m_j$.
\end{Conj}
\end{sloppypar}

An equivalent generating function formulation of
Conjecture~\ref{conj:CyclicInsertion} follows.
\begin{Conj}\label{conj:CyclicInsertionGF}
Let $x_0,x_1,\dots$ be a sequence of commuting indeterminates.
Then
\[
\begin{split}
   &\sum_{n=0}^\infty y^{2n}
   \sum_{\substack{m_j\ge 0\\0\le j\le 2n}}
   \mathscr{C}(m_0,m_1,\dots,m_{2n})
   \ShuffProd{j=0}{2n} x_j^{m_j}\\
   &\qquad = \sum_{m=0}^\infty Z(m)
   \sum_{n=0}^{\lfloor m/2\rfloor} y^{2n}
   (x_0+x_1+\cdots+x_{2n})^{m-2n}.
\end{split}
\]
\end{Conj}
To see the equivalence of Conjectures~\ref{conj:CyclicInsertion}
and~\ref{conj:CyclicInsertionGF}, observe that by the multinomial
theorem,
\begin{eqnarray*}
   &&\sum_{m=0}^\infty Z(m)\sum_{n=0}^{\lfloor m/2\rfloor}
       y^{2n}(x_0+x_1+\cdots+x_{2n})^{m-2n}\\
   &=&\sum_{n=0}^\infty y^{2n}\sum_{m\ge 2n} Z(m)
       (x_0+x_1+\cdots+x_{2n})^{m-2n}\\
   &=&\sum_{n=0}^\infty y^{2n}\sum_{m\ge 2n}Z(m)
      \sum_{m_0+\cdots+m_{2n}=m-2n}\binom{m_0+\cdots+m_{2n}}
      {m_0,\dots,m_{2n}} \prod_{j=0}^{2n} x_j^{m_j}\\
   &=&\sum_{n=0}^\infty y^{2n}\sum_{m\ge 2n} Z(m)
      \sum_{m_0+\cdots+m_{2n}=m-2n}
      \ShuffProd{j=0}{2n} x_j^{m_j}.
\end{eqnarray*}
Now compare coefficients. Although
Conjecture~\ref{conj:CyclicInsertion} remains unproved, it is
nevertheless possible to reduce the problem to that of
establishing the invariance of ${\mathscr{C}}(\vec s)$ for $\vec
s\in C_{2n+1}(m-2n)$.  More specifically, we have the following
non-trivial result.
\begin{theorem}\label{thm:T}  For all non-negative integers $m$ and
$n$ with $m\ge 2n$,
\[
   \sum_{\vec s\in C_{2n+1}(m-2n)}{\mathscr{C}}(\vec s) =
   Z(m)\times|C_{2n+1}(m-2n)| = Z(m)\binom{m}{2n}.
\]
\end{theorem}

\begin{demo} {Example 4} If $m=2n$, Theorem~\ref{thm:T} states that
\[
   {\mathscr{C}}(\{0\}^{2n+1}) = (2n+1)\z(\{3,1\}^n) = Z(2n),
\]
which is equivalent to the Broadhurst-Zagier formula~(\ref{Z31})
(Theorem 1 of~\cite{BBBLc}).
\end{demo}

\begin{demo} {Example 5} If $m=2n+1$, Theorem~\ref{thm:T} states that
\[
   (2n+1){\mathscr{C}}(1,\{0\}^{2n}) = (2n+1)Z(2n+1),
\]
which is Theorem 2 of~\cite{BBBLc}.
\end{demo}

For $m>2n+1$, Theorem~\ref{thm:T} gives new results, although no
additional instances of Conjecture~\ref{conj:CyclicInsertion} are
settled.  For the record, we note the following restatement of
Theorem~\ref{thm:T} in terms of $Z$-functions:
\begin{corollary}[Equivalent to Theorem~\ref{thm:T}]\label{cor:T}
Let $T_{m,n}$ be as in Definition~\ref{def:T}, and put $a=dt/t$,
$b=dt/(1-t)$.   Then, for all non-negative integers $m$ and $n$,
with $m\ge 2n$,
\begin{equation}
   \sum_{\vec s\in C_{2n+1}(m-2n)}Z(\vec s)
   =\iu T_{m,n} = \frac{Z(m)}{2n+1}\binom{m}{2n}
   = \frac{2\pi^{2m}}{(2m+2)!}\binom{m+1}{2n+1}.
\label{T-Integral}
\end{equation}
\end{corollary}

\begin{demo} {Proof of Theorem~\ref{thm:T}}
In view of the equivalent
reformulation~(\ref{T-Integral}) and the well-known
evaluation~(\ref{z2repeated}) for $Z(m)$, it suffices to prove
that with $T_{m,n}$ as in Definition~\ref{def:T} and with
$a=dt/t$, $b=dt/(1-t)$, we have
\[
   \iu T_{m,n} =  \frac{2\pi^{2m}}{(2m+2)!}\binom{m+1}{2n+1}.
\]
Let
\[
   J(z) := \sum_{k=0}^\infty z^{2k} \iu (ab)^k =
   \sum_{k=0}^\infty z^{2k}\z(\{2\}^k).
\]
Then~\cite{BBB} $J(z) = (\sinh(\pi z))/(\pi z)$ for $z\ne 0$ and
$J(0)=1$.  We have
\begin{eqnarray}
   J(z\cos\theta)J(z\sin\theta)
   &=& \frac{\sinh(\pi z\cos\theta)}{\pi z\cos\theta}
   \cdot\frac{\sinh(\pi z\sin\theta)}{\pi z\sin\theta}\nonumber\\
   &=& \frac{\cosh\pi z(\cos\theta+\sin\theta)
   -\cosh\pi z(\cos\theta-\sin\theta)}{2\pi^2 z^2
   \sin\theta\cos\theta}\nonumber\\
   &=& \frac{\cosh\pi z\sqrt{1+\sin 2\theta}-\cosh\pi
   z\sqrt{1-\sin 2\theta}}{\pi^2 z^2 \sin 2\theta}\nonumber\\
   &=&\sum_{m=1}^\infty
   \frac{(\pi z)^{2m}\left\{(1+\sin 2\theta)^m-(1-\sin
   2\theta)^m\right\}}{(2m)!\,\pi^2 z^2 \sin 2\theta}\nonumber\\
   &=&\sum_{m=0}^\infty \frac{2(\pi z)^{2m}}{(2m+2)!}
   \sum_{n=0}^{\lfloor m/2\rfloor} \binom{m+1}{2n+1}(\sin 2\theta)^{2n}.
\label{thm:T-RHS}
\end{eqnarray}
On the other hand, putting $x=z^2\cos^2\theta$ and
$y=z^2\sin^2\theta$ in Theorem~\ref{thm:T-Binom} yields
\begin{eqnarray}
   &&J(z\cos\theta)J(z\sin\theta)\nonumber\\
   &=& \bigg(\sum_{k=0}^\infty (z\cos\theta)^{2k}\iu (ab)^k\bigg)
   \bigg(\sum_{j=0}^\infty (z\sin\theta)^{2j}\iu (ab)^j\bigg)\nonumber\\
   &=& \sum_{m=0}^\infty \sum_{n=0}^m (z\cos\theta)^{2n}
   (z\sin\theta)^{2m-2n} \iu (ab)^n \shuff (ab)^{m-n}\nonumber\\
   &=& \sum_{m=0}^\infty \sum_{n=0}^{\lfloor m/2\rfloor}
   (4 z^4 \sin^2 \theta \cos^2 \theta)^n
   (z^2\cos^2\theta+z^2\sin^2\theta)^{m-2n}\iu T_{m,n}\nonumber\\
   &=& \sum_{m=0}^\infty z^{2m}\sum_{n=0}^{\lfloor m/2\rfloor}
   (\sin 2\theta)^{2n} \iu T_{m,n}.
\label{thm:T-LHS}
\end{eqnarray}
Equating coefficients of $z^{2m}(\sin 2\theta)^{2n}$
in~(\ref{thm:T-RHS}) and~(\ref{thm:T-LHS}) completes the proof.
$\qed$
\end{demo}







\begin{acknowledgment}
Thanks are due to the referee whose comments helped improve the
exposition.
\end{acknowledgment}


\end{article}
\end{document}